\documentclass[a4paper,14pt]{article}

\usepackage{pstricks}
\usepackage{amsmath}
\usepackage{amsfonts}
\usepackage{pstcol}
\usepackage{amscd}

\newcommand{\B}[1]{\mathbb #1}
\newcommand{\C}[1]{{\cal #1}}

\newcommand{\al}{{\alpha}}

\newcommand{\be}{{\beta}}

\newcommand{\Om}{{\Omega}}
\newcommand{\om}{{\omega}}
\newcommand{\eps}{{\varepsilon}}

\newcommand{\Ga}{{\Gamma}}

\newcommand{\La}{{\Lambda}}
\newcommand{\si}{{\sigma}}

\newcommand{\map}[1]{\stackrel {#1}\longrightarrow}

\newcommand{\qed}{\rightline {$\Box $}}

\newcommand{\Mo}{(M,\omega )}

\newcommand{\Ham}{Ham(M,\omega )}

\newcommand{\Mham}{M_{Ham(M,\omega )}}

\newcommand{\BS}{{\bigskip}}
\newcommand{\NI}{{\noindent}}

\newtheorem{theorem}{Theorem}[section]
\newtheorem{thm}[theorem]{Theorem}
\newtheorem{cor}[theorem]{Corollary}

\newtheorem{rem}[theorem]{Remark}
\newtheorem{lemma}[theorem]{Lemma}
\newtheorem{prop}[theorem]{Proposition}
\newtheorem{ex}[theorem]{Example}
\numberwithin{equation}{section}
\numberwithin{figure}{section}

\title{Evaluation fibrations and topology of symplectomorphisms}

\author{Jaros\l aw K\c edra\\
        University of Szczecin
\thanks{
The author is a member of EDGE, Research Training Network
HPRN-CT-2000-00101, supported by The European Human Potential Programme.\newline
Keywords: rational homotopy; symplectic manifold.\newline
AMS classification(2000): Primary 55P62;  Secondary 57R17}}

\date{\today}

\begin{document}


\maketitle

\begin{abstract}
There are two main results. The first states that isotropy
subgroups of groups acting transitively on a rationally 
hyperbolic spaces have infinitely generated rational
cohomology algebra. 
Using this fact, we prove that the analogous statement holds for
groups of symplectomorphisms
of certain blow-ups.
\end{abstract}

\BS

\section{Introduction}\label{S:intro}

Let $G$ be a topological group acting effectively and
transitively on a topological space $X$.  Denote by
${ G}_{pt}$ the isotropy subgroup for the point $pt\in X$.
Let's consider the  following 
{\bf evaluation fibration} associated with the action

$${ G} _{pt} \to { G} \map{ev} X.$$

\noindent
Here $ev(\phi ) = \phi (pt)$, where $pt \in X$ is a point.
Our first aim is to prove the following

\begin{thm}\label{T:iso}
If a topological group $ G$ acts transitively on a simply connected
rationally hyperbolic space $X$, then the rational cohomology
of the isotropy subgroup of a point $H^*( G_{pt};\B Q)$
is infinitely generated as an algebra.
\end{thm}

The proof uses several known facts from
rational homotopy theory and here we sketch the argument.
The first step is to show that the map induced
by the evaluation on homotopy groups
has a finite dimensional image. This is
achieved by the application of the
Gottlieb theory (Theorem \ref{T:got})
under the assumption that the space $X$ has finite
Lusternik-Schnirelmann category.
The second step is to use
the existence of an exact sequence of, so called,
dual rational homotopy groups (Theorem \ref{T:les}) associated
to the fibration. Then it follows from the first step
that all but the finite dimensional part of the rational
homotopy of $X$ comes from the rational homotopy
of the isotropy subgroup.
Since $X$ is rationally hyperbolic,
i.e. its rational homotopy is infinitely dimensional,
then so does the rational homotopy of the isotropy subgroup.
Topological groups are
H-spaces, thus their dual rational homotopy generates the
rational cohomology and we obtain that the latter must
be infinitely generated as an algebra.

\BS

As an application of the above result we prove the following

\begin{thm}\label{T:infty}
Let $\Mo $ be a compact simply connected
symplectic 4-manifold which is neither
rational nor ruled surface up to blow up. Suppose that
$b_2=\dim H^2(M) > 2$. Then the rational cohomology 
$H^*(Symp(\widetilde {M_{\eps }}))$ of the symplectomorphism
group of the symplectic blow up
of $\Mo $ is infinitely generated (as an algebra),
for every sufficiently small
$\eps > 0$.
\end{thm}

\BS
\NI
In all known (very specific) examples, the rational
cohomology rings of symplectomorphism groups are
finitely generated \cite{ab,am}.  
The above result shows that  topology
of groups of symplectomorphisms is complicated in general.

\BS
\NI
{\bf Acknowledgments:} 
The paper was written during my stay at MPI in Bonn and IMPAN in Warsaw.
The final version was completed in the University of Munich. 
I would like to thank Dusa McDuff for very influential discussions.
Conversations with Denis Auroux, Kai Cieliebak, Dieter Kotschick
and Fran\c cois Lalonde 
were also very helpful.

\section{Rational homotopy methods}\label{S:rat}

This section is devoted to recalling  some facts and notions
of rational homotopy theory. We refer to \cite{ap,fe,fht}
for detailed exposition.

\subsection{Minimal models}\label{SS:mm}

\noindent
Given connected topological space X, there is associated a free differential 
graded algebra $\C M_X =(\Lambda V,d)$ such that

\begin{enumerate}

\item there exists a basis $\{v_{\alpha}\}_{\alpha\in J} $

of V for some well ordered set $J$ such that

$$ d(v_\alpha)\in \Lambda V_{<\alpha} $$
and
$$deg(v_\alpha)<deg(v_{\beta})\Rightarrow \alpha<\beta, $$
where $ \Lambda V_{<\alpha} $ is a subalgebra in $ \Lambda V $generated by all
$v_\beta $ for $\beta <\alpha $,

\item there exists a DGA-morphism $ \varphi _X:\C M_X\to \ A(X)$ which induces an isomorphism
on cohomology. Here $A(X)$ denotes Sullivan's rational forms.
$\C M_X$ is called a {\bf minimal  model} of $X$.
\end{enumerate}

The idea of minimal models can be extended to  fibrations.
Namely, one may try to write a minimal model of the total space
of a fibration in terms of the minimal models of the fiber and the
base. While this fails in general, however, it gives a certain model
of the fibration as stated in the following result [AP, Theorem (2.5.1)].

\begin{thm} [Grivel-Halperin-Thomas] \label{T:fght}
         Let $\pi :E\to B$ be a Serre fibration of path connected spaces
         and $F = \pi ^{-1}(b)$ be the fiber over the base-point $b$.
         Suppose that:

       \begin{enumerate}
           \item $F$ is path connected,
           \item $\pi _1(B)$
           acts nilpotently on $H^k(F)$ for all $k\geq 1$,
           \item either $B$ or $F$ has finite $\B K$-type
\footnote{
          Recall that path connected $X$ has finite $\B K$-type if
          $H^k(X;\B K)$ is finite dimensional for  $k\geq 1$.}.
          \end{enumerate}

        \noindent
        Then there exists a KS-model
        \footnote{KS stands for Koszul and Sullivan}
        of the fibration
        \bigskip
$$ 
\CD
           A(B) @>A(\pi )>>   A(E)   @>A(i)>>    A(F) \\
            @AA\varphi _BA          @AA\varphi A       @AA\varphi _F A \\
       {\C M}_B @>>> \C M_B\otimes \C M_F @>>>    \C M_F
\endCD
$$

      \bigskip
      \noindent
      in which $\varphi _B:\C M_B \to A(B)$ is a minimal model of $B$,
      $\varphi _F :\C M_F\to A(F)$ is a minimal model of $F$ and
       $\varphi:\C M_B\otimes \C M_F\to A(E)$ induces an isomorphism
on cohomology.
  \end{thm}

\qed

\BS
\noindent
\subsection{Dual rational homotopy groups}\label{SS:drh}

Let $\C M_X = (\La V, d)$ be a minimal model of the space
$X$. Recall that $V = \bigoplus _{k>0} V^k $ is the graded rational
vector space. The subspaces $V^k$ are called the {\bf dual
rational homotopy groups} of $X$ and denoted by $\Pi ^k(X)$.
Notice that $\Pi ^*$ is a contravariant functor from
the category of path connected topological spaces to the
category of graded rational vector spaces (see \cite{ap}
for details). The next two theorems are also taken from the book
\cite{ap} (Theorem (2.3.7) and Corollary (2.5.2)).

\begin{thm}\label{T:pi}
Let $X$ be a nilpotent space of finite $\B Q$-type 
\footnote{
A path connected space $X$ is said to be nilpotent if $\pi _1(X)$
is a nilpotent group and acts nilpotently on higher homotopy
groups. It has finite $\B Q$-type if $H^n(X;\B Q)$ are finite
dimensional for all $n\geq 1$.}.
Then for each $k\geq 2$, there is a natural isomorphism

$$\Pi^k(X):=Hom_{\B Z}(\pi _k(X),\B Q).$$
Furthermore, this holds for $k=1$, provided that $\pi _1(X)$ is
abelian.\qed
\end{thm}

\begin{thm}\label{T:les}
Under the conditions of Theorem \ref{T:fght}, there is a long exact sequence
of dual homotopy groups

$$
\CD
0\to \Pi^1(B) \to \Pi^1(E) \to \Pi^1(F) \to \Pi^2(B)\to \dots \\
\to \Pi^k(B) \to \Pi^k(E) \to \Pi^k(F) \to \Pi^{k+1}(B)\to \dots 
\endCD
$$
\qed
\end{thm}

\noindent
\subsection{Gottlieb theory}\label{SS:gt}

Let $X$ be a path connected topological space and $HE(X)$
be the space of self-homotopy equivalences of $X$.
Then we have  the evaluation
map $ev:HE(X)\to X$. The images of the homomorphisms
induced by the evaluation on homotopy groups were first
studied by Gottlieb \cite{got1,got2}, now
called the Gottlieb subgroups. We denote them
by $G_k(X) := im(\pi _k(ev))\subset \pi_k(X)$.

To formulate the next result we need the notion of
the Lusternik-Schnirelmann category of a space.
The {\bf LS category} of $X$, denoted $cat(X)$, is the
least integer $m$ (or $\infty $), such that $X$
is the union of $m+1$ open subsets contractible in $X$.
The following theorem is an easy consequence of 
Proposition 28.8 in \cite{fht}.

\begin{thm}\label{T:got}
Suppose that $X$ is a simply connected topological space
of finite category. Then 
\begin{enumerate}
\item $G_{2k} \otimes \B Q = 0$ and
\item dim$G_* \otimes \B Q \leq cat(X)$.
\end{enumerate}\qed
\end{thm}

\section{The rational dichotomy and proof of Theorem \ref{T:iso}}\label{S:res}

Suppose that $X$ is simply connected
and has the homotopy type of a finite CW-complex.
Let $G$ is a connected group
acting transitively on $X$ with the isotropy
group of a point $pt\in X$ denoted by $G_{pt}$.

\begin{thm}\label{T:main}
With the above notation,
there is the following estimate
 
$$\dim \pi _k(X)\otimes \B Q \leq
  \dim \pi _k(G_{pt})\otimes \B Q,$$
for $k$ big enough.
\end{thm}

\noindent
{\bf Proof:} The assumption are made in order to apply
Theorem \ref{T:fght} and \ref{T:les} to the evaluation
fibration

$$
\CD
G_{pt} \to G \to X.
\endCD
$$

\noindent
Thus Theorem \ref{T:les} implies that there exists
the following exact sequence 
$$
\dots \to \Pi^{k-1}(G_{pt}) \to \Pi^ k(X) \to \Pi ^k(G)
\to \Pi ^k(G_{pt}) \to \dots
$$

\noindent
Since $X$ is up to homotopy a finite CW-complex, then it
has finite LS-category. It follows from Theorem \ref{T:got}
that $\dim G_*(X)\oplus \B Q $ is finite, hence
$ev^*:\Pi ^k(X)\to \Pi ^k(G)$ is trivial for $k$ big enough.
Now the statement follows from the exactness of the
sequence of dual homotopy groups. 

\qed

\BS

The rational dichotomy discovered by F\'elix \cite{fe,fht}
states that if $X$ is an n-dimensional simply connected
space with the rational cohomology of finite type
and finite category, then either its rational homotopy
is finite dimensional  or the dimensions of rational
homotopy groups
grow exponentially. In the first case the space is
called {\bf rationally elliptic} and in the second 
{\bf rationally hyperbolic}. As in other branches of
mathematics, hyperbolicity is in a sense a "generic" feature.
This is illustrated in the following

\begin{prop}\label{P:ell}
If $X$ is rationally elliptic then
\begin{enumerate}
\item $\dim \pi _{even}\otimes \B Q \leq 
       \dim \pi _{odd}\otimes \B Q \leq cat(X)$,
\item $\chi(X)\geq 0$, where $\chi $ denotes the 
       Euler characteristic.\qed
\end{enumerate}
\end{prop}

\noindent
{\bf Proof of Theorem \ref{T:iso}:} 
Since $G_{pt}$ is a topological group
then its rational cohomology is freely generated its by dual 
rational homotopy. If $X$ is rationally hyperbolic then
Theorem \ref{T:main} implies that the dual rational
homotopy is of infinite dimension.

\qed

\BS

\begin{ex}\label{E:4dim}
{\em
Let $X$ be a 4-dimensional simply connected finite CW-complex.
Then it follows form the basic properties of the 
Lusternik-Schnirelmann category that $cat(X)\leq 2$. 
Now the above proposition implies
that if $\pi _2(X)\otimes \B Q \geq 3$, then $X$ is
rationally hyperbolic. 
Thus every simply connected 4-manifold $X$, whose 
$b^2(X):=\dim H^2(X;\B Q) \geq 3$ is rationally hyperbolic.
By Theorem \ref{T:iso}, we get that if a topological group
acts  transitively on $X$ then the isotropy subgroup has
infinitely generated rational cohomology.}
\end{ex}

\section {The cohomology of a symplectomorphism group is
infinitely generated}\label{S:infty}

\subsection {The argument}

In this section we prove Theorem \ref{T:infty}. We need to introduce
some additional notions. Let $B_{\eps }$ be the standard symplectic
ball of capacity $\eps $ and $Symp\Mo$ denote the component of
the identity in the group of symplectomorphisms of $\Mo $.
Notice that since $M$ is simply connected, then 
$Symp\Mo = Ham\Mo $ the group of Hamiltonian symplectomorphisms.
We also need the following:

\begin{itemize}


\item $Symp^{U(2)}(M,B_{\eps })$ is the subgroup of $Symp\Mo $
consisting of elements acting $U(2)$-linearly on $B_{\eps}$.

\item $Emb^{U(2)}(B_{\eps},M)$ is the set of symplectic embeddings
of the $\eps$-ball into $\Mo $ modulo the $U(2)$ reparametrizations
of the source.
\end{itemize}

\NI
Consider the following commutative diagram of fibrations

$$
\CD
Symp^{U(2)}(M,B_{\eps}) @>>> Symp(M,pt) \\
  @VVV                       @VVV       \\
Symp\Mo        @=            Symp\Mo    \\
 @VVV                    @V{ev}VV       \\
Emb^{U(2)}(B_{\eps},M)  @>{f}>> M          
\endCD
$$


\BS
\NI
The argument is split into several lemmas, which we shall prove
in the next section.

\begin{lemma}\label{L:J-curves}
Let $\Mo $ be a compact symplectic manifold which is
neither  rational nor ruled surface up to blow up.
Let $(\widetilde{M_{\eps }},\om_{\eps})$ be the symplectic
blow up of $\Mo $.
Then for
any almost complex structure $J$ compatible with $\om_{\eps} $
there exist unique $J$-holomorphic exceptional sphere 
which is embedded.
\end{lemma}

\begin{lemma}\label{L:blow}
If for
any almost complex structure $J$ compatible with $\om_{\eps} $
there exist unique $J$-holomorphic exceptional sphere 
which is embedded, then 
$Symp(\widetilde{M_{\eps}})$.
is weak homotopy equivalent
to $Symp^{U(2)}(M,B_{\eps})$. 
\end{lemma}


\begin{lemma}\label{L:fsur}
In the above diagram, the map $f$ induces a surjection on
rational homotopy groups, except $\pi_2$ on which it
depends on the first Chern class of $\Mo$.
\end{lemma}

\BS
\NI
{\bf Proof of Theorem \ref{T:infty}:}
Consider the long exact sequences of homotopy groups.

$$
\CD
\pi_{k+1}(Emb^{U(2)}(B_{\eps},M))  @>{f_*}>>\pi_{k+1}( M)\\
          @VVV                       @VVV       \\          
\pi_k(Symp^{U(2)}(M,B_{\eps})) @>>>\pi_k( Symp(M,pt)) \\
          @VVV                       @VVV       \\
\pi_k(Symp\Mo) @=                  \pi_k(Symp\Mo )   \\
@VVV                                @V{ev_*}VV       \\
\pi_k(Emb^{U(2)}(B_{\eps},M))  @>{f_*}>>\pi_k( M)          
\endCD
$$

\BS
\NI
Theorem \ref{T:got} implies that 
the rank of the image of the connecting homomorphism
$\pi_{*+1}(M)\to \pi_*(Symp(M,pt))$ is infinite. Thus,
according to Lemma \ref{L:fsur}, the same is true for
$\pi_{*+1}(Emb^{U(2)}(B_{\eps},M))\to \pi_*(Symp^{U(2)}(M,B_{\eps}))$.
In particular,
$\dim \pi_*(Symp^{U(2)}(M,B_{\eps}))\otimes \B Q = \infty$.
The statement of the theorem follows immediately from
Lemma \ref{L:blow} and the proof of Theorem \ref{T:iso}.

\qed


\subsection {Proofs of the lemmas}

\BS
\NI
{\bf Proof of Lemma \ref{L:J-curves}:} This essentially Theorem 4.39
in \cite{ms}. First, it is clear that there exists exceptional sphere
which is embedded for any generic $J$. According to the criterion
of Hofer, Lizan and Sikorav \cite{hls} we see that it is regular.
It is also unique due to positivity of intersections.

Suppose that for some $J$ the class $E$ of the exceptional sphere
is represented by some cusp-curve $\sum C_i$. Then 
$1=C_1(E) = \sum c_1(C_i)$ which implies that for some $i$
$c_1(C_i)\geq 0$. Now Theorem 4.39 in \cite{ms} says that $\Mo $ is
a blow up of a rational or ruled surface which is excluded
by the hypothesis  (cf. Proposition 2.6 in \cite{lp}).

\qed

\BS
\NI
{\bf Proof of Lemma \ref{L:blow}:} This was proved by Lalonde
and Pinsonnault in \cite{lp}  Lemma 2.3 and 2.4.

\qed

\BS
\NI
We need some preparation to prove Lemma \ref{L:fsur}.
Let $J$ be an almost complex structure on $M$ compatible with $\om$.
Let $B_{\eps}\to TM_{\eps} \to M$ be the $\eps$-ball subbundle of the
tangent bundle (with respect to the Riemannian metric given by $J$).
We need the following parametrized version of the Darboux theorem:

\begin{prop}[Universal thickening property]\label{P:univ}
There exist a map $\phi : TM_{\eps}\to M$ such 
that for
any $p\in M$ the composition 
$\phi \circ i_m:B_{\eps}\to M$
is a symplectic embedding.
Here $i_m:B_{\eps}\to TM_{\eps}$
is an inclusion of the fiber over the point $m\in M$.
\end{prop}

\NI
{\bf Proof:}
We start with the exponential map
associated to the Riemannian metric,

$$
\exp:TM_{\eps}\to M.
$$

\NI
Then $\exp ^*\om \in \Om^2(TM_{\eps })$ and we have a map

$$
V:\Om^2(TM_{\eps})\to \Om^2_{vert}(TM_{\eps}),
$$
where the latter space consist of 2-forms evaluated only
on the vectors tangent to the fibers of $TM_{\eps}$.
In other words, this is a space of sections
of the vector bundle over $M$ whose
fibers are 2-forms on the fibers of $TM_{\eps}$, that is

$$
\Om^2_{vert}(TM_{\eps}) = \Ga (M,P\times _{U(n)} \Om^2(B_{\eps})),
$$
where $P\to M$ is $U(n)$-principal frame bundle of $M$.
The differential $d_B$ on $B_{\eps}$ defines a morphism
of vector bundles

$$
d:P\times _{U(n)} \Om^1(B_{\eps}) \to P\times _{U(n)} \Om^2(B_{\eps}) 
$$
by 
$d[p,\al ] = [p, d_B(\al )]$. It is well defined because
$U(n)$ acts linearly on $B_{\eps }$. Since the ball $B_{\eps}$
is contractible then the map $d$ is surjective.

The bundle $P\times _{U(n)} \Om^2(B_{\eps})$ has an obvious
section $\widehat {\om_0} :M\ni m\mapsto [p, \om_0]$, where
$\om_0$ is the standard symplectic form on $B_{\eps}$.
Let's consider the section

$$
V(\exp^*\om) - \widehat {\om _0} \in \Ga (M,P\times _{U(n)} \Om^2(B_{\eps})).
$$
We can assume (after possibly composing exp with a linear automorphism
of $TM_{\eps}$) that the above section is equal to zero when restricted
to the zero section of $TM_{\eps}$.

Because the map $d$ is surjective we can find a section $\si $ of
$P\times _{U(n)} \Om^1(B_{\eps})$ such that 
$d\circ \si = V(\exp^*\om) - \widehat {\om_0}$. Moreover,
$\si $ can be chosen so that it vanishes along the zero
section of $TM_{\eps}$. Indeed, if $\si '$ is any such section,
then we take $\si := \si ' - \si_0'$, where $\si_0' $ over $m\in M$
is the constant 1-form 
equal to $\si ' (m)$ at the origin. It is clearly closed, so
$d\circ \si = 
d\circ (\si ' - \si_0' ) = V(\exp^*\om) - \widehat {\om_0} = a'$. 

Now, in the pesence of symplectic form in each fiber
(the constant section $\om_0$) the section $\si $
gives rise to the section of

$$
P\times _{U(n)}\mathfrak{X} (B_{\eps}),
$$
where $\mathfrak{X} (B_{\eps})$ denotes vector fields on $B_{\eps}$.
In other words we get a vector field on $TM_{\eps}$ tangent
to the fibers. Moreover, this vector field is trivial along
the zero section. Take the flow $\psi^t$ of this vector field
(and smaller $\eps $ if necessary). We obtain that

$$
V(\exp \circ \psi ^1)^*\om = \widehat {\om_0}.
$$

\NI
Notice that when restricted to a fiber of $TM_{\eps}$
the above agrument reduces to the standard proof of the
Darboux theorem \cite[Theorem 3.15]{ms}.
Finally we define $\phi := \exp \circ \, \psi^1$. Then for
an inclusion of the fiber $i_m:B_{\eps}\to TM_{\eps}$ we
get that $\exp \circ \, \psi^1 \circ i_m:B_{\eps}\to M$
is a Darboux chart, that is a symplectic embedding.

\qed

\BS
In general we define
a {\bf twisted family of 
embeddings of $F$ into $M$ parametrized by $B$}
to be a map
$\phi :E\to M$,
where $F\to E\to B$ is a bundle and $\phi $ restricted
to every fiber is an embedding.
This is in contrast to the usual notion of a family of
maps where the domain is always a product.
Thus the map 
$\phi: TM_{\eps}\to M$ 
 in the above proposition is an example of a
twisted family of {\em symplectic} embeddings.
This notion is proved useful below in the proof of Lemma \ref{L:fsur}.
The idea is that after restricion of the parameter space
the twisted family becomes trivial, i.e. a usual family.

\begin{cor}[Thickening property]\label{C:univ}
Let $f:X\to \Mo $ be a continuous (or smooth) map. Then there exist
a continuous (or smooth) twisted family of symplectic embeddings 
of small balls $B_{\eps}$
parametrized by $X$. 
\end{cor}

\NI
{\bf Proof:}
Take the pull back bundle 

$$
\CD
f^*TM_{\eps} @>{\widehat{f}}>> TM_{\eps}\\
@VVV                             @VVV   \\
X            @>f>>             M.
\endCD
$$

\NI
The twisted family of symplectic embeddings is now defined by
$\phi \circ \widehat {f}:f^*TM_{\eps} \to M$.

\qed

\begin{rem}
{\em
We call the above facts {\em thickening property} because
they say that given any map we can always ``thicken'' it
to get a twisted family of embeddings of symplectic balls.
}
\end{rem}

\BS
\NI
{\bf Proof of Lemma \ref{L:fsur}:}
Let $s\in \pi_k(M)$. According to Corollary \ref{C:univ}, we have
a map $\phi :s^*TM_{\eps} \to M$ such that its restriction to
each fiber is symplectic embedding. If the bundle $s^*TM_{\eps }$
is trivial then it defines an element of $\pi_k(Emb(B_{\eps},M))$.
Notice that the rational homotopy of the structure group of
$s^*TM_{\eps}$ (which is $Sp(4,\B R)$) is equal to
the exterior algebra
$\Lambda (c_1,c_2)$, where $c_i$ is of degree $2i-1$.
Thus in degrees different from 2 and 4 we get that
$ks^*(TM_{\eps})$ is trivial, where $k\in \B Z$.
It means that 

$$
f_*:\pi_k(Emb(B_{\eps},M))\otimes \B Q\to \pi_k(M)\otimes \B Q 
$$
is surjective for $k\neq 2,4$. In degree 4 it is also
true because there is no nonzero degree map
$S^4\to M$, where $M$ is symplectic 4-dimensional manifold. 
In  degree 2, we have an obvious
dependence on the corresponding Chern classes.

\qed

\subsection {Final remarks}

\begin{enumerate}
\item It is very likely that the assumption in Theorem \ref{T:infty} 
saying that the manifold is not blow up of neither rational or ruled
surface can be weakened. Then one has to put more work to show that
for every compatible almost complex structure there exist embedded
exceptional curve (cf. Proposition 2.6 in \cite{lp}).

\item
Theorem \ref{T:infty} can be proved in slightly another way.
Namely, it is possible to prove that the map
$\pi_k(Symp^{U(2)}(M,B_{\eps})) \to \pi_k( Symp(M,pt))$
is a surjection. This together with Lemma \ref{L:blow}
proves the theorem. I learned this argument from
Fran\c cois Lalonde.

\item Notice that the infinite dimensional part
of the rational homotopy of the symplectomorphism group
detected by Theorem \ref{T:infty} is robust. That is
it survives when we embed symplectomorphisms into
homeomorphisms or even into homotopy equivalences.

\item Although the rational cohomology ring of the group
of symplectomorphisms is infinitely generated, it does not
mean that the topology of this group cannot be understood.
The hope is that there is more structure. Namely, there
is the Pontriagin product on homology and with respect
to it the homology may be finitely generated. 
Thus the more appropriate structure to investigate 
in topology of symplectomorphism groups is
is their homology with the Pontriagin product. 

\item 
Similarly, cohomology ring of the classifying space
need not be infinitely generated in this case. If it is
the case then it means that there is a lot of nontrivial
fibrations with trivial characteristic classes.

\end{enumerate}

\BS
\BS

Jarek K\c edra

Institute of Mathematics US

Wielkopolska 15

70-451 Szczecin

Poland

\end{document}